\documentclass[11pt,leqno]{amsart}
\usepackage{amsmath,epsfig,graphicx,color}
\definecolor{darkblue}{rgb}{.2, 0.2,.8}
\definecolor{carageen}{rgb}{0,0.5,0.3}
\definecolor{darkred}{rgb}{.8, .1,.1}
\newcommand{\red}{\color{darkred}}

\newcommand{\slln}{strong law of large numbers}

\newcommand{\spr}{stochastic process}

\newcommand{\ex}{{\rm e}\,}

\newcommand{\asy}{asymptotic}

\newtheorem{lemma}{Lemma}[section]

\newtheorem{theorem}[lemma]{Theorem}

\newtheorem{proposition}[lemma]{Proposition}
\newtheorem{definition}[lemma]{Definition}
\newtheorem{corollary}[lemma]{Corollary}
\newtheorem{example}[lemma]{Example}
\newtheorem{exercise}[lemma]{Exercise}
\newtheorem{remark}[lemma]{Remark}
\newtheorem{fig}[lemma]{Figure}
\newtheorem{tab}[lemma]{Table}

\newcommand{\bfR}{{\bf R}}

\newcommand{\bth}{\begin{theorem}}
\newcommand{\ethe}{\end{theorem}}

\newcommand{\bre}{\begin{remark}\em }
\newcommand{\ere}{\end{remark}}

\newcommand{\ble}{\begin{lemma}}
\newcommand{\ele}{\end{lemma}}

\newcommand{\bde}{\begin{definition}}
\newcommand{\ede}{\end{definition}}
\newcommand{\bco}{\begin{corollary}}
\newcommand{\eco}{\end{corollary}}

\newcommand{\bpr}{\begin{proposition}}
\newcommand{\epr}{\end{proposition}}

\newcommand{\bexer}{\begin{exercise}}
\newcommand{\eexer}{\end{exercise}}

\newcommand{\bexam}{\begin{example}}
\newcommand{\eexam}{\end{example}}

\newcommand{\bfi}{\begin{fig}}
\newcommand{\efi}{\end{fig}}

\newcommand{\btab}{\begin{tab}}
\newcommand{\etab}{\end{tab}}

\newcommand{\fidi}{finite-dimensional distribution}
\newcommand{\rv}{random variable}

\newcommand{\cov}{{\rm cov}}

\newcommand{\rhs}{right-hand side}

\newcommand{\beao}{\begin{eqnarray*}}
\newcommand{\eeao}{\end{eqnarray*}\noindent}

\newcommand{\beam}{\begin{eqnarray}}
\newcommand{\eeam}{\end{eqnarray}\noindent}

\newcommand{\beqq}{\begin{equation}}
\newcommand{\eeqq}{\end{equation}\noindent}

\newcommand{\bce}{\begin{center}}
\newcommand{\ece}{\end{center}}

\newcommand{\barr}{\begin{array}}
\newcommand{\earr}{\end{array}}
\newcommand{\cadlag}{c\`adl\`ag}

\newcommand{\stp}{\stackrel{\P}{\rightarrow}}

\newcommand{\stas}{\stackrel{\rm a.s.}{\rightarrow}}

\newcommand{\eqd}{\stackrel{d}{=}}

\newcommand{\vague}{\stackrel{\lower0.2ex\hbox{$\scriptscriptstyle
                    \it{v} $}}{\rightarrow}}
\newcommand{\weak}{\stackrel{\lower0.2ex\hbox{$\scriptscriptstyle
                    \it{w} $}}{\rightarrow}}
\newcommand{\what}{\stackrel{\lower0.2ex\hbox{$\scriptscriptstyle
                    \it{\hat{w}} $}}{\rightarrow}}

\newcommand{\bdis}{\begin{displaymath}}
\newcommand{\edis}{\end{displaymath}\noindent}

\newcommand{\N}{\mathbb{N}}
\newcommand{\R}{\mathbb{R}}

\newcommand{\nto}{n\to\infty}

\newcommand{\xto}{x\to\infty}

\newcommand{\ov}{\overline}
\newcommand{\wt}{\widetilde}

\newcommand{\vep}{\varepsilon}

\newcommand{\regvary}{regularly varying}
\newcommand{\slvary}{slowly varying}

\newcommand{\bfT}{{\bf T}}
\newcommand{\bbr}{{\mathbb R}}

\newcommand{\bbz}{{\mathbb Z}}

\newcommand{\BM}{Brownian motion}

\newcommand{\con}{convergence}

\newcommand{\st}{such that}
\newcommand{\fif}{if and only if}

\newcommand{\chf}{characteristic function}
\newcommand{\fct}{function}

\newcommand{\ds}{distribution}

\newcommand{\rep}{representation}

\newcommand{\seq}{sequence}

\newcommand{\pro}{probabilit}

\newcommand{\ms}{measure}

\newcommand{\ld}{large deviation}
\newcommand{\bfx}{{\bf x}}
\newcommand{\bfX}{{\bf X}}

\newcommand{\bfY}{{\bf Y}}

\newcommand{\bft}{{\bf t}}

\newcommand{\bfs}{{\bf s}}

\newcommand{\E }{{\mathbb E}}
\renewcommand{\P }{{\mathbb P}}

\allowdisplaybreaks 
\begin{document}
\today
\bibliographystyle{plain}
\title[Distance covariance for stochastic processes]{Distance covariance for stochastic processes}
\thanks{
Muneya Matsui's research is partly supported
by JSPS Grant-in-Aid
for Young Scientists B (16K16023) and Nanzan University
Pache Research Subsidy I-A-2 for the 2016 academic year.
Thomas Mikosch's research is partly supported
by the Danish Research Council Grant DFF-4002-00435. 
Gennady Samorodnitsky's research  is partly supported by the ARO MURI
grant W911NF-12-1-0385.}
\author[M. Matsui]{Muneya Matsui}
\address{Department of Business Administration, Nanzan University, 18
Yamazato-cho, Showa-ku, Nagoya 466-8673, Japan.}
\email{mmuneya@gmail.com}
\author[T. Mikosch]{Thomas Mikosch}
\address{Department  of Mathematics,
University of Copenhagen,
Universitetsparken 5,
DK-2100 Copenhagen,
Denmark}
\email{mikosch@math.ku.dk}
\author[G. Samorodnitsky]{Gennady Samorodnitsky}
\address{School of Operations Research and Information Engineering,
Cornell University,
220 Rhodes Hall,
Ithaca, NY 14853, U.S.A.}
\email{gennady@orie.cornell.edu}

\begin{abstract}
The distance covariance of two random vectors is a \ms\ of their dependence.
The empirical distance covariance and correlation can be used as statistical tools for
testing whether two random vectors are independent.
We propose an analogs
of the distance covariance for two stochastic processes defined on some interval. Their empirical
analogs can be used to test the independence of two processes. 
\end{abstract}
\keywords{Empirical \chf , distance covariance, stochastic process, test of independence}
\subjclass{Primary 62E20; Secondary 62G20 62M99 60F05 60F25}
\maketitle

The authors of this paper would like to congratulate Tomasz Rolski
on his 70th birthday. We would like to express our gratitude for his longstanding
contributions to applied probability theory as an author, editor, and organizer.
Tomasz kept applied probability going in Poland and beyond even in difficult historical times. 
The applied probability community, including ourselves, has benefitted a lot 
from his enthusiastic, energetic and reliable work.
\bce
{\bf Sto lat! 
 Niech zyje nam! 
 Zdrowia, szczescia, pomyslnosci!} 
\ece

\section{Distance covariance for processes on $[0,1]$}\label{sec:model}\setcounter{equation}{0}
We consider a real-valued stochastic process $X=(X(t))_{t\in [0,1]}$
with sample paths in a measurable space $S$ such that $X$ is
measurable as a map from its probability space into $S$. We assume
that the  \pro y \ms\    $P_X$ generated by $X$ on $S$ is uniquely
determined by its \fidi s. Examples include processes with continuous
or \cadlag\ sample paths on $[0,1]$. 
The \pro y  measure $P_X$ is then determined by the totality of the
\chf s 
\beao
\varphi_X(\bfx_k;\bfs_k)=\varphi^{(k)}_X(\bfx_k;\bfs_k)
= \int_{S} \ex^{i\,(s_1\,f(x_1)+\cdots + s_k\,f(x_k))}\,P_X(df)\,,\qquad k\ge 1\,,
\eeao
where $\bfx_k=(x_1,\ldots,x_k)'\in [0,1]^k$,
$\bfs_k=(s_1,\ldots,s_k)'\in\bbr^k$. In particular, for two such
processes,  $X$ and $Y$,  the measures $P_X$ and $P_Y$ coincide if and
only if 
\beao
\varphi_X(\bfx_k;\bfs_k)=\varphi_Y(\bfx_k;\bfs_k)\qquad\mbox{for all $\bfx_k\in[0,1]^k,\bfs_k\in\bbr^k$, $k\ge 1$.}
\eeao
\par
We now turn from the general question of identifying the \ds s of $X$ and $Y$ to a more 
specific but related one: given two processes
$X,Y$ on $[0,1]$ with values in $S$ as above and 
defined on the same \pro y space, we intend to find some  means to
verify whether $X$ and $Y$ are independent. Motivated by the
discussion above, we need to show that the joint  law of $(X,Y)$ on
$S\times S$, denoted by $P_{X,Y}$, coincides with the product  \ms\
$P_X\otimes P_Y$. Assuming, once again, that a probability measure on
$S\times S$ is determined by the finite-dimensional distributions (as
is the case with the aforementioned examples), we need to show that the joint
\chf s of $(X,Y)$ factorize, i.e., 
\beam\label{eq:1}
\varphi_{X,Y}(\bfx_k;\bfs_k,\bft_k)&=& \int_{S^2} \ex^{i\,\sum_{j=1}^k (s_j f(x_j)+t_jh(x_j))}\,P_{X,Y}(df,dh)\nonumber\\
&=&\varphi_X(\bfx_k;\bfs_k) \,\varphi_Y(\bfx_k;\bft_k)
\,,\qquad\bfx_k\in [0,1]^k,\bfs_k,\bft_k\in\bbr^k\,,k\ge 1\,.
\eeam
Clearly, this condition is hard to check and therefore we try to get a more compact equivalent condition which can also
be used for some statistical test of independence between $X$ and $Y$.
\par
For this reason, we consider  a unit rate Poisson process
$N=(N(t))_{t\in [0,1]}$ with arrivals $0<T_1<T_2<\cdots <T_{N(1)}\le 1$, 
write $\bfT_N=(T_1,\ldots,T_{N(1)})'$ and, correspondingly
$\bfs_N,\bft_N$ for any vectors in $\bbr^{N(1)}$. Then, for any
positive \pro y density \fct\ $g$ on $\bbr$, we define 
\beam\label{eq:2}
\lefteqn{d(P_{X,Y},P_X\otimes P_Y)}\nonumber\\&=& 
\E_N\Big[\int_{\bbr^{2\,N(1)}}
\big|\varphi_{X,Y}(\bfT_N;\bfs_N,\bft_N)-\varphi_X(\bfT_N;\bfs_N)\,\varphi_Y(\bfT_N;\bft_N)\big|^2\,
\prod_{j=1}^{N(1)}g(s_j) g(t_j)\, d \bfs_N\,d\bft_N\Big]\nonumber\\
&=& \sum_{k=1}^\infty \P(N(1)=k)\,\int_{[0,1]^{k}}\Big[\int_{\bbr^{2k}}
\big|\varphi_{X,Y}(\bfx_k;\bfs_k,\bft_k)-\varphi_X(\bfx_k;\bfs_k)\varphi_Y(\bfx_k;\bft_k)\big|^2\,\nonumber\\
&&\hspace{4cm}\times \prod_{j=1}^kg(s_j)g(t_j)\, d \bfs_kd\bft_k\,\Big]
d\bfx_k\,, 
\eeam
where in the last step we used  the order
statistics property of the homogeneous Poisson  process. 
Here we interpret the summand corresponding to $k=0$ as zero, and we also suppress the dependence on $g$ in the notation. 
Now, the right-hand integrals vanish \fif\ \eqref{eq:1} is satisfied for Lebesgue a.e. $\bfx_k,\bfs_k,\bft_k$, hence \fif\ \eqref{eq:1}
holds for any $\bfx_k,\bfs_k,\bft_k$. We summarize:
\ble\label{lem:1} If $g$ is a positive probability density on $\bbr$ then
$d(P_{X,Y},P_X\otimes P_Y)=0$ \fif\ $P_{X,Y}=P_X\otimes P_Y$.
\ele
\bre
Lemma~\ref{lem:1} can easily be extended in several directions.\\
1. The statement remains valid when
the Poisson \pro ies $(\P(N(1)=k))_{k\ge 1 }$ are replaced by any
summable \seq\ of  nonnegative numbers with infinitely many positive terms.\\
2.  Obvious modifications of Lemma \ref{lem:1} are valid e.g. for random
fields $X,Y$ on $[0,1]^d$ 
(in this case we can sample the values
of the random fields at the points of a Poisson random measure on
$[0,1]^d$ whose mean measure is the $d$-dimensional Lebesgue
measure). Moreover, the values of $X,Y$ may be multivariate. \\
3. The positive \pro y density $\prod_{j=1}^kg(s_j)g(t_j)$ on $\bbr^{2k}$ can be replaced by any positive measurable \fct\ provided
the infinite series in \eqref{eq:2} is finite.
 This idea will be exploited in Section~\ref{sec:4} below.\\  
4. Returning
to our original problem about identifying the laws of $X$ and $Y$, similar calculations 
show that the quantity
\beao
d(P_X,P_Y)&=&
\sum_{k=1}^\infty \P(N(1)=k)\,\int_{[0,1]^{k}}\Big[\int_{\bbr^{k}}
\big|\varphi_{X}(\bfx_k;\bfs_k)-\varphi_Y(\bfx_k;\bfs_k)\big|^2\,\prod_{j=1}^kg(s_j)\, d \bfs_k\Big] d\bfx_k\nonumber\\
\eeao
vanishes  \fif\ $X\eqd Y$, where $\eqd$ means that all
finite-dimensional distributions of $X$ and $Y$ coincide.
The quantity $d(P_X,P_Y)$ can be taken as the basis for
a goodness-of-fit test for the \ds s of $X$ and $Y$. 
\ere
\par
In what follows, we refer to the quantities $d(P_{X,Y},P_X\otimes P_Y)$ as {\em distance covariance} between the stochastic processes
$X$ and $Y$. This name is motivated by work on {\em distance covariance} for random 
vectors $\bfX\in\bbr^p,\bfY\in\bbr^q$ (possibly of different dimensions) defined by 
\beao
T(\bfX,\bfY)= \int_{\bbr^{p+q}} \big|\varphi_{\bfX,\bfY}(\bfs,\bft)-\varphi_{\bfX}(\bfs)\,\varphi_{\bfY}(\bft)\big|^2\,\mu(d\bfs,d\bft)\,,
\eeao
where $\mu$ is a (possibly infinite) \ms\ on $\bbr^{p+q}$; see for example
\cite{dmmw:2016,feuerverger:1993,szekely:rizzo:bakirov:2007,szekely:rizzo:2009,szekely:rizzo:2014}.
The last mentioned authors coined the names {\em distance covariance} and {\em distance correlation} for the
standardized version $R(\bfX,\bfY)=T(\bfX,\bfY)/\sqrt{T(\bfX,\bfX)T(\bfY,\bfY)}$; they chose some special infinite \ms s $\mu$ which 
lead to an elegant form
of $T(\bfX,\bfY)$ and $R(\bfX,\bfY)$; see Section~\ref{sec:4} for more information on this approach. The goal in the aforementioned literature was to find a statistical tool for testing independence
between the vectors $\bfX$ and $\bfY$ using the fact that $R(\bfX,\bfY)=0$ \fif\ $\bfX,\bfY$ are independent provided $\mu$ has a positive
Lebesgue density on $\bbr^{p+q}$. The sample versions $T_n(\bfX,\bfY)$ and $R_n(\bfX,\bfY)=T_n(\bfX,\bfY)/\sqrt{T_n(\bfX,\bfX)T_n(\bfY,\bfY)}$, 
constructed from an iid sample $(\bfX_i,\bfY_i)$, $i=1,\ldots,n,$ of
copies of $(\bfX,\bfY)$, are then used as test statistics for checking independence of $\bfX$ and $\bfY$.
\par
For stochastic processes $X,Y$ on $[0,1]$ one might be tempted to test their independence based on independent 
observations $\bfX_i=(X_i(x_1),\ldots,X_i(x_k))'$, $\bfY_i=(Y_i(x_1),\ldots,Y_i(x_k))'$, $i=1,\ldots,n$ of the processes $X,Y$ at the locations $\bfx_k$ in
$[0,1]^k$. However, \cite{szekely:rizzo:2013} observed that the empirical 
distance correlation $R_n(\bfX,\bfY)$ has the tendency to be very close to 1 even for relatively small values $k$.
Our approach avoids the high dimensionality of the vectors $\bfX_i$ and $\bfY_i$ by randomizing their dimension~$k$.
\par
Our paper is organized as follows. In Section~\ref{sec:2} we study some of the theoretical properties of
the distance covariance between two \spr es $X,Y$ on $[0,1]$ where we assume that $g$ is a positive \pro y density.
We find a tractable \rep\ of this distance covariance from which we derive the corresponding sample version. 
In Section~\ref{sec:4} we choose
the non-integrable weight \fct\ $g$ from \cite{szekely:rizzo:bakirov:2007}. Again,
we find a suitable \rep\ of this distance covariance, derive the corresponding sample version and show
that it is a consistent estimator of its deterministic counterpart. In Section~\ref{sec:5}
we conduct a small simulation study based on the sample distance correlation introduced in Section~\ref{sec:2}.
We compare the small sample behavior of the sample distance correlation with the corresponding sample distance correlation of 
\cite{szekely:rizzo:bakirov:2007} for independent and dependent Brownian and fractional Brownian sample paths.

\section{Properties of distance covariance}\label{sec:2}\setcounter{equation}{0} 
\subsection{Distance correlation}
In the context of \spr es $X,Y$ one may be interested in standardizing the distance covariance $T(X,Y)=d(P_{X,Y},P_X\otimes P_Y)$, i.e.,
in the {\em distance correlation}
\beao
R(X,Y)= \dfrac{T(X,Y)}{\sqrt{T(X,X)\,T(Y,Y)}}\,.
\eeao
However, it is not obvious that $R(X,Y)$   assumes only values between
$0$ and $1$. This property is guaranteed  by a Cauchy-Schwarz
argument.  
\ble
Assume  that $g(s)=g(-s)$. Then $0\leq R(X,Y)\leq 1$.
\ele
We have $R(X,X)=1$, In general, the relation $R(X,Y)=1$ does not imply $X=Y$ a.s. For example, if $X$ is symmetric 
then $R(X,-X)=1$ as well.
\begin{proof}
 Let $(X^\prime, Y^\prime)$ be an independent copy of $(X,Y)$.
Applying   the Cauchy-Schwarz inequality first to the $k$-dimensional
integral with respect to the 
product of $k$ copies of $g$, then to the
expectation with respect to the law of $(X,Y)$, next with respect to
the Lebesgue measure on $[0,1]^k$ and, finally, with respect to the
law of $N$, and using the symmetry of the density $g$, we obtain 
\beao
T(X,Y) &=&
\sum_{k=1}^\infty \P(N(1)=k)\,\int_{[0,1]^{k}} d\bfx_k \\
&& \E\Big[\int_{\bbr^{2k}} 
\left[\left( \ex^{i\sum_{j=1}^k s_jX_j} -  \varphi_X(\bfx_k;\bfs_k)\right) 
\left( \ex^{i\sum_{j=1}^k t_jY_j} -  \varphi_Y(\bfx_k;\bft_k)\right)
\right.\\
 && \left.  
\left( \ex^{-i\sum_{j=1}^k s_jX^\prime_j} -  \varphi_X(\bfx_k;-\bfs_k)\right) 
\left( \ex^{-i\sum_{j=1}^k t_jY^\prime_j} -  \varphi_Y(\bfx_k;-\bft_k)\right)
\right] \prod_{j=1}^kg(s_j)g(t_j)\, d \bfs_kd\bft_k\Big] \\
&\leq& \sum_{k=1}^\infty \P(N(1)=k)\,\int_{[0,1]^{k}} d\bfx_k \\
&& \Big( \E \Big[\Big|\int_{\bbr^{k}}  \Big( \ex^{i\sum_{j=1}^k s_jX_j} -  \varphi_X(\bfx_k;\bfs_k)\Big) 
\Big( \ex^{-i\sum_{j=1}^k s_jX^\prime_j} -   \varphi_X(\bfx_k;-\bfs_k)\Big)  
\prod_{j=1}^kg(s_j)\, d \bfs_k \Big|^2\Big]\Big)^{1/2}  \\
&&   \times\ \Big(\E\Big[\Big| \int_{\bbr^{k}}  \Big( \ex^{i\sum_{j=1}^k t_jY_j} -  \varphi_Y(\bfx_k;\bft_k)\Big) 
\Big( \ex^{-i\sum_{j=1}^k t_jY^\prime_j} -   \varphi_Y(\bfx_k;-\bft_k)\Big)  
\prod_{j=1}^kg(t_j)\, d \bft_k \Big|^2\Big]\Big)^{1/2}  \\
&=& \sum_{k=1}^\infty \P(N(1)=k)\,\int_{[0,1]^{k}} d\bfx_k \\
&&\Big[\int_{\bbr^{2k}}
\big|\varphi_{X,X}(\bfx_k;\bfs_k,\bft_k)-\varphi_X(\bfx_k;\bfs_k)\varphi_X(\bfx_k;\bft_k)\big|^2\,\prod_{j=1}^kg(s_j)g(t_j)\, d \bfs_kd\bft_k\,\Big]^{1/2}\\&&
\times \Big[\int_{\bbr^{2k}}
\big|\varphi_{Y,Y}(\bfx_k;\bfs_k,\bft_k)-\varphi_Y(\bfx_k;\bfs_k)\varphi_Y(\bfx_k;\bft_k)\big|^2\,\prod_{j=1}^kg(s_j)g(t_j)\, d \bfs_kd\bft_k\,\Big]^{1/2}\\
&\le &\sum_{k=1}^\infty \P(N(1)=k)\,\\&&\Big[\int_{[0,1]^{k}} d\bfx_k\int_{\bbr^{2k}}
\big|\varphi_{X,X}(\bfx_k;\bfs_k,\bft_k)-\varphi_X(\bfx_k;\bfs_k)\varphi_X(\bfx_k;\bft_k)\big|^2\,\prod_{j=1}^kg(s_j)g(t_j)\, d \bfs_kd\bft_k\,\Big]^{1/2}\\&&
\times \Big[\int_{[0,1]^{k}} d\bfx_k\int_{\bbr^{2k}}
\big|\varphi_{Y,Y}(\bfx_k;\bfs_k,\bft_k)-\varphi_Y(\bfx_k;\bfs_k)\varphi_Y(\bfx_k;\bft_k)\big|^2\,\prod_{j=1}^kg(s_j)g(t_j)\, d \bfs_kd\bft_k\,\Big]^{1/2}\\
&\le &  \sqrt{T(X,X)}\,\sqrt{T(Y,Y)}\,.
\eeao
This proves that $0\leq R(X,Y)\leq 1$. 
\end{proof}
\subsection{Representations}
Our next goal is to find explicit expressions for $d(P_{X,Y},P_X\otimes P_Y)$.
We observe that
\beao
\lefteqn{\big|\varphi_{X,Y}(\bfx_k;\bfs_k,\bft_k)-\varphi_X(\bfx_k;\bfs_k)\,\varphi_Y(\bfx_k;\bft_k)\big|^2}\\
&=& |\varphi_{X,Y}(\bfx_k;\bfs_k,\bft_k)|^2 + |\varphi_X(\bfx_k;\bfs_k)|^2|\varphi_Y(\bfx_k;\bft_k)|^2 -
2\,{\rm Re}\,\{\varphi_{X,Y}(\bfx_k;\bfs_k,\bft_k)\varphi_X(\bfx_k;-\bfs_k)\varphi_Y(\bfx_k;-\bft_k)\}\,.
\eeao
This expression suggests to decompose \eqref{eq:2} into 3 distinct
parts, the first one being 
\beao
&&\sum_{k=1}^\infty \dfrac{\ex^{-1}}{k!}\,\int_{[0,1]^k}\Big[\int_{\bbr^{2k}} |\varphi_{X,Y}(\bfx_k;\bfs_k,\bft_k)|^2
\,\prod_{j=1}^kg(s_j)g(t_j)\, d \bfs_k\,d\bft_k\Big] d\bfx_k\\
&=& \int_{S^2}\sum_{k=1}^\infty \dfrac{\ex^{-1}}{k!}\Big(\int_{[0,1]^{k}}\Big[ \int_{\bbr^{2k}}
\ex^{i\sum_{r=1}^k \big(s_r \,(f(x_r)-f'(x_r))+ t_r \,(h(x_r)-h'(x_r))\big)}\\
&&\hspace{2cm}\prod_{j=1}^k g(s_j) g(t_j) \,d\bfs_k d\bft_k \Big]d\bfx_k\Big) P_{X,Y}(d(f,h))\,P_{X,Y}(d(f',h'))\\
&=& \int_{S^2} \sum_{k=1}^\infty \dfrac{\ex^{-1}}{k!} \Big(
\int_{[0,1]} \Big[\int_{\bbr} 
\ex^{is \,(f(x)-f'(x))}\,g(s) ds\,\int_{\bbr} 
\ex^{it \,(h(x)-h'(x))}\,g(t) dt\Big] dx\Big)^k\\
&&\hspace{2cm}P_{X,Y}(d(f,h))\,P_{X,Y}(d(f',h'))\\
&=&\ex^{-1}\int_{S^2}\Big[ \exp\Big(\int_{[0,1]} 
\Big[\int_{\bbr^2} \ex^{is \,(f(x)-f'(x))+ it\,(h(x)-h'(x))} 
g(s)\,g(t)\, 
ds\,dt\Big] dx\Big)-1  \Big]\\&&
\hspace{2cm}P_{X,Y}(d(f,h))\,P_{X,Y}(d(f',h'))\,.
\eeao
Similar calculations yield
\beao
d(P_{X,Y},P_X\otimes P_Y)&=&\ex^{-1}\int_{S^2}\Big[\exp\big(\int_{[0,1]} \int_{\bbr} 
\ex^{is \,(f(x)-f'(x))} g(s)\,ds\,\int_{\bbr}
\ex^{is\, (h(x)-h'(x))} g(s)\,ds dx \big)\Big]\\
&&\times\,
\big[P_{X,Y}(d(f,h))\,P_{X,Y}(d(f',h'))+ P_X\otimes P_Y(d(f,h)) P_X\otimes P_Y(d(f',h'))\\
&&-P_{X,Y}(d(f,h))\,P_X\otimes P_Y(d(f',h'))-P_{X,Y}(d(f',h'))\,P_{X}\otimes P_Y (d(f,h))\big]\,.
\eeao
We summarize our results:
\ble\label{lem:m}
The distance covariance between the processes $X,Y$ on $[0,1]$ with values in $S$ can be written
in the following form:
\beao
\ex^{1}d(P_{X,Y},P_X\otimes P_Y)&=&
\E \Big[\exp\big(\int_{[0,1]} \int_\bbr 
\ex^{is \,(X(x)-X'(x))}\,g(s)\,ds\int_\bbr \ex^{is (Y(x)-Y'(x))} g(ds)ds\, dx \big)\Big]\\
&&+\E \Big[\exp\big(\int_{[0,1]} \int_{\bbr} 
\ex^{is \,(X(x)-X'(x))}g(s) ds \int_\bbr\ex^{is (Y''(x)-Y'''(x))}g(s)\,ds  dx \big)\Big]\\
&&-2{\rm Re}\E \Big[\exp\big(\int_{[0,1]} \int_{\bbr} 
\ex^{is \,(X(x)-X'(x))}g(s) ds \int_\bbr\ex^{is (Y(x)-Y''(x))}g(s)\,ds dx \big)\Big]\,,
\eeao
where $(X',Y')$ is an independent copy of $(X,Y)$ and $Y'',Y'''$ are independent copies of $Y$ which are also independent of $X,X',Y,Y'$.
\ele
\bexam \rm Let $g$ be the density of  a suitably scaled symmetric $\alpha$-stable law on $\bbr$, $\alpha\in (0,2]$. Then 
\beao
\int_{\bbr}\ex^{is \,(f(x)-f'(x))} g(s)\,ds =
\ex^{-|f(x)-f'(x)|^\alpha} \,, 
\eeao
and so for a uniform \rv\ $U$ on $(0,1)$ which is independent of $X,Y,X',Y',Y'',Y'''$, 
\beam
d(P_{X,Y},P_X\otimes P_Y)&=&\ex^{-1}
\E \Big[\exp\big(\E_U \ex^{-|X(U)-X'(U)|^\alpha-|Y(U)-Y'(U)|^\alpha }\big)\nonumber\\
&&+\exp\big(\E _U \ex^{-|X(U)-X'(U)|^\alpha-|Y''(U)-Y'''(U)|^\alpha}\big)\nonumber\\
&&-2\,\exp\big(\E_U \ex^{-|X(U)-X'(U)|^\alpha-|Y(U)-Y''(U)|^\alpha}\big)\Big]\,,\label{eq:kk}
\eeam
where $E_U$ denotes expectation with respect to $U$.
\eexam
\subsection{Sample distance covariance}
Let $(X_1,Y_1),\ldots,(X_n,Y_n)$ be an iid sample with \ds\ $P_{X,Y}$ and let $P_{n,X,Y}$ be the corresponding empirical
\ds\ with marginals $P_{n,X}$ and $P_{n,Y}$. Then we can define the {\em sample distance covariance} given by
\beao
T_n(X,Y)&=&\ex^1\,d(P_{n,X,Y},P_{n,X}\otimes P_{n,Y} )\\
&=& \dfrac 1 {n^2} \sum_{j_1=1}^n\sum_{j_2=1}^n 
\exp\Big(\int_{[0,1]} \int_\bbr \ex^{is \,(X_{j_1}(x)-X_{j_2}(x))} g(s)\,ds 
\int_\bbr \ex^{is \,(Y_{j_1}(x)-Y_{j_2}(x))} g(s)\,ds 
 dx \Big)\\
&&+\dfrac 1 {n^4}\sum_{j_1=1}^n\sum_{j_2=1}^n  \sum_{j_3=1}^n\sum_{j_4=1}^n  
\exp\Big(\int_{[0,1]} \int_\bbr \ex^{is \,(X_{j_1}(x)-X_{j_2}(x))} g(s)
\,ds  \int_\bbr \ex^{is \,(Y_{j_3}(x)-Y_{j_4}(x))} g(s) ds
dx \Big)\\
&&-2\,{\rm Re}\dfrac 1 {n^3}\sum_{j_1=1}^n\sum_{j_2=1}^n \sum_{j_3=1}^n  \exp\Big(\int_{[0,1]} \int_\bbr \ex^{is \,(X_{j_1}(x)-X_{j_2}(x))}g(s)\,ds \int_\bbr \ex^{is \,(Y_{j_1}(x)-Y_{j_3}(x))}g(s)\,ds 
dx \Big)\,.
\eeao
\bre
This estimator is the exact sample analog of the distance covariance.
However, this estimator is of $V$-statistics-type and  leads to an additional 
bias.
For practical purposes, one should avoid summation over diagonal and subdiagonal terms, making
the estimator of $U$-statistics-type. Then, for example, the first expression would turn into
\beao
\dfrac 1 {n(n-1)} \sum_{j_1=1}^n\sum_{j_2=1,j_2\ne j_1}^n   
\exp\Big(\int_{[0,1]} \int_\bbr \ex^{is \,(X_{j_1}(x)-X_{j_2}(x))} g(s)\,ds 
\int_\bbr \ex^{is \,(Y_{j_1}(x)-Y_{j_2}(x))} g(s)\,ds 
 dx \Big)\,.
\eeao
Since  the bias is \asy ally negligible and 
we are interested only in \asy\ results we stick
to the original version of the sample distance covariance. In Section~\ref{sec:4} we consider a distinct 
version of distance covariance; see \eqref{eq:distinct}. By virtue of its construction diagonal and subdiagonal terms
vanish in its sample version, i.e., a bias problem does not appear.
\ere

\bexam Assume that $g$ is the density of a suitably scaled symmetric $\alpha$-stable \rv . Then
\beao
\lefteqn{\ex^1\,d(P_{n,X,Y},P_{n,X}\otimes P_{n,Y} )}\\
&=&\dfrac 1 {n^2} \sum_{j_1=1}^n\sum_{j_2=1}^n  \exp\Big(\int_{[0,1]} \ex^{-|X_{j_1}(x)-X_{j_2}(x)|^\alpha-|Y_{j_1}(x)-Y_{j_2}(x)|^\alpha}\,dx\Big)\\
&&+\dfrac 1 {n^4}\sum_{j_1=1}^n\sum_{j_2=1}^n  \sum_{j_3=1}^n\sum_{j_4=1}^n  
\exp\Big(\int_{[0,1]}  \ex^{-|X_{j_1}(x)-X_{j_2}(x)|^\alpha-|Y_{j_3}(x)-Y_{j_4}(x)|^\alpha}\,
dx\Big)\\
&&-\dfrac 2 {n^3}\sum_{j_1=1}^n\sum_{j_2=1}^n \sum_{j_3=1}^n  \exp\Big(\int_{[0,1]}  
\ex^{-|X_{j_1}(x)-X_{j_2}(x)|^\alpha-|Y_{j_1}(x)-Y_{j_3}(x)|^\alpha}\,dx\Big)\,.
\eeao
\eexam
\bre
The form of the sample distance covariance indicates that one needs to involve numerical methods for its calculation.
In addition, in general we cannot assume that the sample paths of $(X_i,Y_i)$ are completely observed. Then we need to approximate
the path-dependent integrals appearing in the exponents of the
expressions above by appropriate sums on a grid. These problems are
not studied further in this paper.
\ere
\par
The following result is an immediate con\seq\ of the \slln\ for $U$-statistics (see \cite{hoffmann:1994}) and the observation that
$d(P_{n,X,Y},P_{n,X}\otimes P_{n,Y} )$ is a linear combination of $U$-statistics. 
\bpr
Assume that $\big((X_i,Y_i)\big)_{i=1,\ldots,n }$ is an iid \seq\ of $S^2$-valued random elements. Then 
\beao
d(P_{n,X,Y},P_{n,X}\otimes P_{n,Y} )\stas d(P_{X,Y},P_{X}\otimes P_{Y} )\,,\quad \nto\,.
\eeao
\epr
\section{Distance covariance with infinite weight \ms s}\label{sec:4}\setcounter{equation}{0} 
So far we assumed that $g$ is a positive integrable density. In the aforementioned literature (see 
for example \cite{szekely:rizzo:bakirov:2007}) positive weight \fct s $g$ were used which are
not integrable over $\bbr$.
In what follows, we consider an approach with suitable positive non-integrable weight \fct s which lead to a distance covariance  for
\spr es.  Due to positivity of this weight \fct\ Lemma~\ref{lem:1} remains valid.

To begin, note that if the function $g$  is not necessarily integrable
but is symmetric, then 
appealing to \eqref{eq:2} and using the symmetry of both the cosine
\fct\ and the function $g$ we have 
\beam\label{eq:ex}
\lefteqn{d(P_{X,Y},P_X\otimes P_Y)}\nonumber\\
&=& \sum_{k=1}^\infty \P(N(1)=k)\,\int_{[0,1]^{k}}\E\Big[\int_{\bbr^{2k}}\Big(
\cos(\bfs_k '(\bfX_k-\bfX'_k))\cos(\bft_k '(\bfY_k-\bfY'_k))\nonumber\\&&
+\cos(\bfs_k '(\bfX_k-\bfX'_k))\cos(\bft_k '(\bfY_k''-\bfY'''_k))
- 2\cos(\bfs_k '(\bfX_k-\bfX'_k))\cos(\bft_k '(\bfY_k-\bfY''_k)) \Big)\nonumber\\
&& \hskip 2in\prod_{j=1}^kg(s_j)g(t_j)\, d \bfs_kd\bft_k\,\Big] d\bfx_k\,,
\eeam
where $\bfX_k= (X(x_1),\ldots,X(x_k))'$, $\bfY_k= (Y(x_1),\ldots,Y(x_k))'$ and $(\bfX_k',\bfY_k')$ is an independent copy of 
 $(\bfX_k,\bfY_k)$ while $\bfY_k'',\bfY_k'''$ 
are iid copies of $\bfY_k$ independent of everything else.
Since 
\beam\label{eq:ex1}
\cos u \cos v =1 -(1-\cos u) -(1-\cos v) + (1-\cos u)(1-\cos v)\,,
\eeam
we have
\beao
\lefteqn{d(P_{X,Y},P_X\otimes P_Y)}\\
&=& \sum_{k=1}^\infty \P(N(1)=k)\,\int_{[0,1]^{k}}\E\Big[\int_{\bbr^{2k}}\Big(
(1-\cos(\bfs_k '(\bfX_k-\bfX'_k)))(1-\cos(\bft_k '(\bfY_k-\bfY'_k)))\\&&
+(1-\cos(\bfs_k '(\bfX_k-\bfX'_k)))(1-\cos(\bft_k '(\bfY_k''-\bfY'''_k)))\\&&
-2 (1-\cos(\bfs_k '(\bfX_k-\bfX'_k)))(1-\cos(\bft_k '(\bfY_k-\bfY''_k)))\Bigl)\\
&& \hskip 2in \prod_{j=1}^kg(s_j)g(t_j)\, d \bfs_kd\bft_k\,\Big] d\bfx_k\,.
\eeao
\par
 Next we replace the product kernels 
$\prod_{j=1}^kg(s_j)$ above by other positive measurable
functions on $\bbr^k$.  Inspired by \cite{szekely:rizzo:bakirov:2007}
we choose the \fct s 
\beao
g_k(\bfs)= c_k\,|\bfs|^{-k-\alpha}\,,\qquad \bfs\in \bbr^k\,,\ \alpha\in (0,2)\,,
\eeao 
where the constant $c_k=c_k(\alpha)>0$ is \st
\beao
\int_{\bbr^k} (1-\cos(\bfs'\bfx))\,g_k(\bfs)\,d\bfs= |\bfx|^\alpha\,,\qquad x\in\bbr^k\,.
\eeao
The corresponding distance covariance between $X$ and $Y$ becomes:
\beao
\lefteqn{d(P_{X,Y},P_X\otimes P_Y)}\\
&=& \sum_{k=1}^\infty \P(N(1)=k)\,\int_{[0,1]^{k}}\E \Big[\int_{\bbr^{2k}}\Big(
(1-\cos(\bfs_k '(\bfX_k-\bfX'_k)))(1-\cos(\bft_k '(\bfY_k-\bfY'_k)))\\&&\quad
+(1-\cos(\bfs_k '(\bfX_k-\bfX'_k)))(1-\cos(\bft_k '(\bfY_k''-\bfY'''_k)))\\&&\quad
- 2(1-\cos(\bfs_k '(\bfX_k-\bfX'_k)))(1-\cos(\bft_k '(\bfY_k-\bfY''_k)))\Bigr)\\
&&\qquad \times g_k(\bfs_k)g_k(\bft_k)\, d \bfs_kd\bft_k\,\Big] d\bfx_k\,.
\eeao
By Fubini's theorem and the order statistics property of the Poisson process,
\beam\label{eq:distinct}
\lefteqn{d(P_{X,Y},P_X\otimes P_Y)}\nonumber\\
&=& \sum_{k=1}^\infty \P(N(1)=k)\,\int_{[0,1]^{k}}
\Big(\E[ |\bfX_k-\bfX'_k|^\alpha |\bfY_k-\bfY'_k|^\alpha]
+\E[|\bfX_k-\bfX'_k|^\alpha]\,\E[|\bfY_k''-\bfY'''_k|^\alpha]\nonumber\\&&
- 2\E[|\bfX_k-\bfX'_k|^\alpha|\bfY_k-\bfY''_k|^\alpha]
\Big) d\bfx_k\nonumber\\
&=& \E [|\bfX_{N}-\bfX'_{N}|^\alpha|\bfY_{N}-\bfY'_{N}|^\alpha]
+\E [|\bfX_{N}-\bfX'_{N}|^\alpha\,|\bfY''_{N}-\bfY'''_{N}|^\alpha]\nonumber\\
&&- 2\,\E[|\bfX_{N}-\bfX'_{N}|^\alpha|\bfY_{N}-\bfY''_{N}|^\alpha]\nonumber\\
&=&I_1+I_2-2\,I_3\,,\nonumber\\
\eeam
where  $\bfX_{N}=(X(T_1),\ldots,X(T_{N(1)}))'$, $\bfY_{N}=(Y(T_1),\ldots,Y(T_{N(1)}))'$,  etc. In particular, all the
expectations are  finite if  
\begin{equation} \label{e:cond}
\sup_{0\leq x\leq 1} \E
[|X(x)|^\alpha+|Y(x)|^\alpha+ |X(x)Y(x)|^\alpha]  <\infty\,.
 \end{equation}
An empirical version of $I_1$ is then given by
\beao
\hat I_1= \dfrac 1 {l_n} \dfrac 1 {n^2}\sum_{1\le i,j\le n} \sum_{k=1}^{l_n} |\bfX_{i,{N_k}}-\bfX_{j,N_k}|^\alpha
|\bfY_{i,N_k}-\bfY_{j,N_k}|^\alpha\,,
\eeao 
where $((X_k,Y_k))$ are iid copies of $(X,Y)$ independent of the iid copies $(N_i)$ of the homogeneous Poisson process $N$. The empirical versions $\hat I_2,\hat I_3$ of $I_2,I_3$ are 
defined in an analogous way. The integer \seq\ $(l_n)$ is \st\ $l_n\to\infty$. 
\par
In view of the \slln\ for $U$-statistics, for fixed $l$, as $\nto$,
\beao\lefteqn{
\dfrac 1 l\sum_{k=1}^l A_{nk}=\dfrac 1 l \dfrac 1 {n^2}\sum_{1\le i,j\le n}\sum_{k=1}^l |\bfX_{i,N_k}-\bfX_{j,N_k}|^\alpha
|\bfY_{i,N_k}-\bfY_{j,N_k}|^\alpha}\\
&\stas& \dfrac 1 l \sum_{k=1}^l\E [|\bfX_{N_k}-\bfX'_{N_k}|^\alpha|\bfY_{N_k}-\bfY'_{N_k}|^\alpha\mid N_k]:= \dfrac 1l\sum_{k=1}^l A_k\,.
\eeao 
Therefore, we can choose a \seq\ $\epsilon_n\downarrow 0$ \st
\beao
\P\Big(\dfrac 1 l \Big|\sum_{k=1}^l (A_{nk}-A_k)\Big|>\epsilon_n\Big)\to 0
\eeao
and then also choose an integer \seq\ $(r_n)$ \st\ $r_n\to\infty$ and 
\beao
r_n\,\P\Big(\dfrac 1 l \Big|\sum_{k=1}^l (A_{nk}-A_k)\Big|>\epsilon_n\Big)\to 0\,.
\eeao
Note that the sequence $(r_n)$ can be chosen to be monotone and such
that $r_n-r_{n-1}\in \{ 0,1\}$ for each $n$. Then 
\beao
\P\Big(\dfrac {1}{r_n\,l} \Big|\sum_{s=1}^{r_n} \sum_{k=(s-1)l+1}^{sl} (A_{nk}-A_k)\Big|>\epsilon_n\Big)
&\le &\P\Big(\dfrac 1 l\sup_{s=1,\ldots,r_n}\Big|\sum_{k=(s-1)l+1}^{sl} (A_{nk}-A_k)\Big|>\epsilon_n\Big)\to 0\,.
\eeao
This means that 
\beao
\dfrac {1}{r_n\,l} \sum_{k=1}^{r_n\,l} (A_{nk}-A_k)\stp 0\,,\qquad \nto\,.
\eeao
However, by the \slln , as $n\to\infty$,
\beao
\dfrac {1}{r_n\,l}  
\sum_{k=1}^{r_n\,l} A_k\stas \E [A_1]=\E [|\bfX_{N}-\bfX'_{N}|^\alpha|\bfY_{N}-\bfY'_{N}|^\alpha]\,.
\eeao
Hence, for every $l$ there is an $(r_n)$ \st
\beao
\dfrac {1}{r_n\,l} \sum_{k=1}^{r_n\,l} A_{nk} \stp \E [A_1]\,,\qquad \nto\,.
\eeao
We conclude that 
\beao
\lefteqn{\sup_{l\,r_{n-1}\le v\le l\,r_n} \Big|\dfrac 1v \sum_{k=1}^v A_{nk}- \dfrac 1{l\,r_n} \sum_{k=1}^{lr_n}A_{nk}\Big|}\\
&\le &\dfrac{r_n-r_{n-1}}{l\,r_{n-1}r_n}  \sum_{k=1}^{l\,r_n} A_{nk}
+ \dfrac 1 {lr_n}\sum_{k=lr_{n-1}+1}^{lr_n} A_{nk}\,.
\eeao
The \rhs\ converges in \pro y to zero, hence we have the law of large numbers for $\hat I_1$. Similar arguments apply to 
$\hat I_2,\hat I_3$. We summarize:
\bpr
Let $\alpha\in (0,2)$ and assume that \eqref{e:cond} holds. 
Then for any integer \seq\ $(l_n)$ with $l_n\to\infty$,
\beao
\lefteqn{d(P_{n,X,Y},P_{n,X}\otimes P_{n,Y})
=\dfrac 1 {l_n} \dfrac 1 {n^2}\sum_{1\le i,j\le n} \sum_{k=1}^{l_n} |\bfX_{i,N_k}-\bfX_{j,N_k}|^\alpha
|\bfY_{i,N_k}-\bfY_{j,N_k}|^\alpha}\\&& +\dfrac 1 {l_n} \dfrac 1 {n^2}\sum_{1\le i,j\le n} \sum_{k=1}^{l_n} |\bfX_{i,N_k}-
\bfX_{j,N_k}|^\alpha \dfrac 1 {l_n} \dfrac 1 {n^2}\sum_{1\le i,j\le n} \sum_{k=1}^{l_n} |\bfY_{i,N_k}-\bfY_{j,N_k}|^\alpha\\
&&-2\dfrac 1 {l_n} \dfrac 1 {n^3}\sum_{1\le i,j,l\le n} \sum_{k=1}^{l_n} |\bfX_{i,N_k}-\bfX_{j,N_k}|^\alpha
|\bfY_{i,N_k}-\bfY_{l,N_k}|^\alpha\\
&\stp & d(P_{X,Y},P_X\otimes P_{Y})\,.
\eeao
\epr

\section{A simulation study}\label{sec:5}\setcounter{equation}{0} 
In what follows, we conduct a small simulation study
for the sample distance correlation $R_n(X,Y)$ from Section~\ref{sec:2} for the standard normal density $g$. This choice implies that
\beao
T_n(X,Y)&=&\ex^1\,d(P_{n,X,Y},P_{n,X}\otimes P_{n,Y} )\\
&=&\dfrac 1 {n^2} \sum_{j_1=1}^n\sum_{j_2=1}^n  \exp\Big(\int_{[0,1]} 
\ex^{-|X_{j_1}(x)-X_{j_2}(x)|^2/2-|Y_{j_1}(x)-Y_{j_2}(x)|^2/2}\,dx\Big)\\
&&+\dfrac 1 {n^4}\sum_{j_1=1}^n\sum_{j_2=1}^n  \sum_{j_3=1}^n\sum_{j_4=1}^n  
\exp\Big(\int_{[0,1]}  \ex^{-|X_{j_1}(x)-X_{j_2}(x)|^2/2-|Y_{j_3}(x)-Y_{j_4}(x)|^2/2}\,
dx\Big)\\
&&-\dfrac 2 {n^3}\sum_{j_1=1}^n \sum_{j_2=1}^n \sum_{j_3=1}^n  \exp\Big(\int_{[0,1]}  
\ex^{-|X_{j_1}(x)-X_{j_2}(x)|^2/2-|Y_{j_1}(x)-Y_{j_3}(x)|^2/2}\,dx\Big)\,.
\eeao
As a matter of fact, simulations of this quantity are highly complex.
We choose a moderate sample size $n=100$ and approximate the integrals on $[0,1]$ by their Riemann sums at an equidistant 
grid with mesh $1/50$. For $(X,Y)$, we take a bivariate Brownian motion $(B_{1},B_{2})$ with correlation $\rho\in[0,1]$, i.e.,  
\[
 \cov(B_1(s),B_2(t)) = \rho \min_{s,t\in[0,1]} \{s,t\}\,,\qquad s,t\in [0,1]\,,
\] 
and a bivariate fractional Brownian motion $(W_1,W_2)$ with correlation
$\rho\in[0,1]$, i.e.,
\[
 \cov(W_1(s),W_2(t)) = \frac{\rho}{2}\{|s|^{2H} + |t|^{2H}
 -|t-s|^{2H}\},\qquad s,t\in[0,1]\,,
\]
where we assume that the Hurst parameters of $W_1$ and $W_2$ are the same; see \cite{lavancier:philippe:surgailis:2009} 
for more general cross-correlation structures of vector-fractional \BM s.
\par
We compare the behavior of the sample distance correlation 
\beao
 R_n(X,Y) =\dfrac{T_n(X,Y)}{\sqrt{T_n(X,X)}\sqrt{T_n(Y,Y)}} 
\eeao
 of the aforementioned \spr es
with the corresponding sample distance correlation from 
 \cite{szekely:rizzo:bakirov:2007}
\beao
R_n^{\rm Sz}(\bfX,\bfY) =\dfrac{T_n^{\rm Sz}(\bfX,\bfY)}{\sqrt{T_n^{\rm Sz}(\bfX,\bfX)}\sqrt{T_n^{\rm Sz}(\bfY,\bfY)}}\,,
\eeao
 where for a sample $(\bfX_i,\bfY_i),i=1,\ldots,n,$ of independent copies of the vector $(\bfX,\bfY)$,
\beao
T_n^{\rm Sz}(\bfX,\bfY)&=& \dfrac 1 {n^2} \sum_{j_1=1}^n\sum_{j_2=1} ^n |\bfX_{j_1}-\bfX_{j_2}||\bfY_{j_1}-\bfY_{j_2}|
+\dfrac{1}{n^4} \sum_{j_1=1}^n\sum_{j_2=1}^n  \sum_{j_3=1}^n\sum_{j_4=1}^n  |\bfX_{j_1}-\bfX_{j_2}||\bfY_{j_3}-\bfY_{j_4}|\\
&&-2 \dfrac 1 {n^3}\sum_{j_1=1}^n\sum_{j_2=1}^n  \sum_{j_3=1}^n  |\bfX_{j_1}-\bfX_{j_2}||\bfY_{j_1}-\bfY_{j_3}|\,.
\eeao
We calculate the sample distance correlation $R_n^{\rm Sz}(\bfX,\bfY)$
based on $n=100$ iid simulations of the vector $(\bfX,\bfY)=(X(i/50),Y(i/50))_{i=1,\ldots,50}$. The calculation of $R_n(X,Y)$ and $R_n^{\rm Sz}(\bfX,\bfY)$
is based on the same simulated sample paths
$((X_i,Y_i))_{i=1,\ldots,n}$. 
\par
\begin{figure}[htbp] 
\begin{center}
\includegraphics[width=1\textwidth,height=10cm]{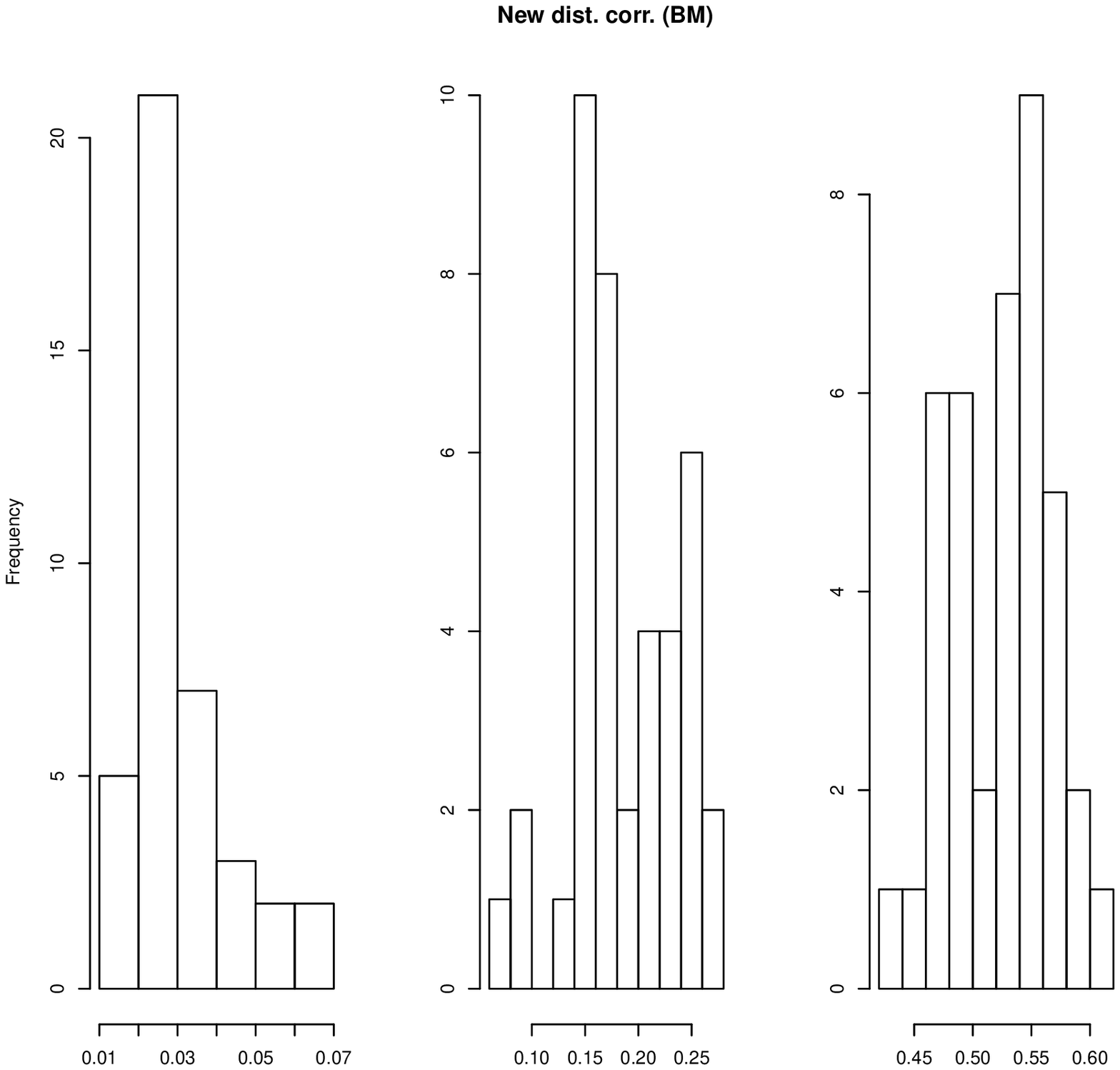}
\includegraphics[width=1\textwidth,height=10cm]{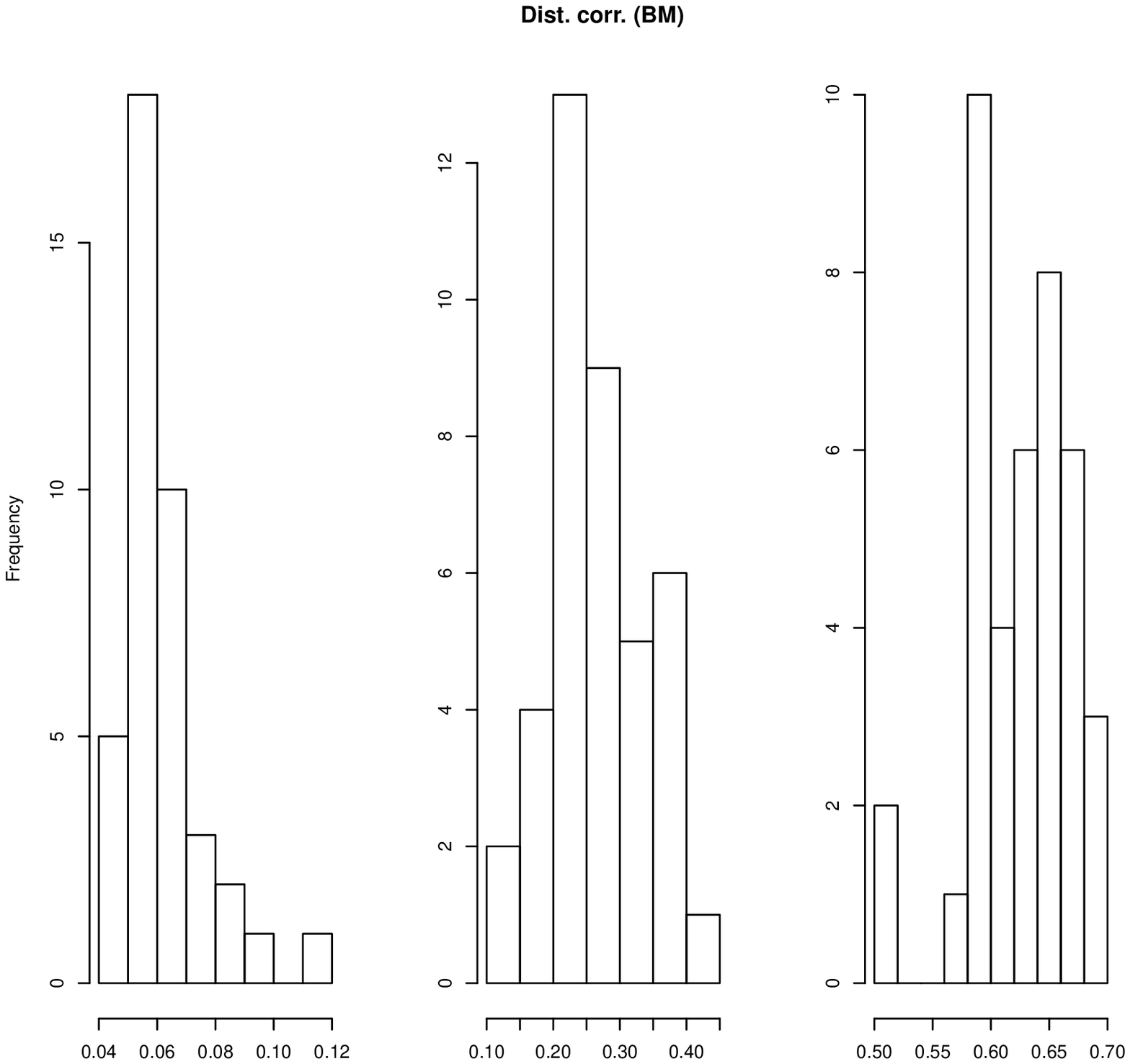} 
\caption{Histograms of $R_n(B_1,B_2)$ (top) and $R_n^{\rm Sz}(B_1,B_2)$ (bottom) based on $40$ samples.  The 
 correlations of $B_1$ and $B_2$ are respectively $\rho=0,\,0.5,\,0.8$,
 from left to right.}
\label{fig:1}
\end{center}
\end{figure}
\begin{figure}[htbp]
\begin{center}
\includegraphics[width=1\textwidth,height=10cm]{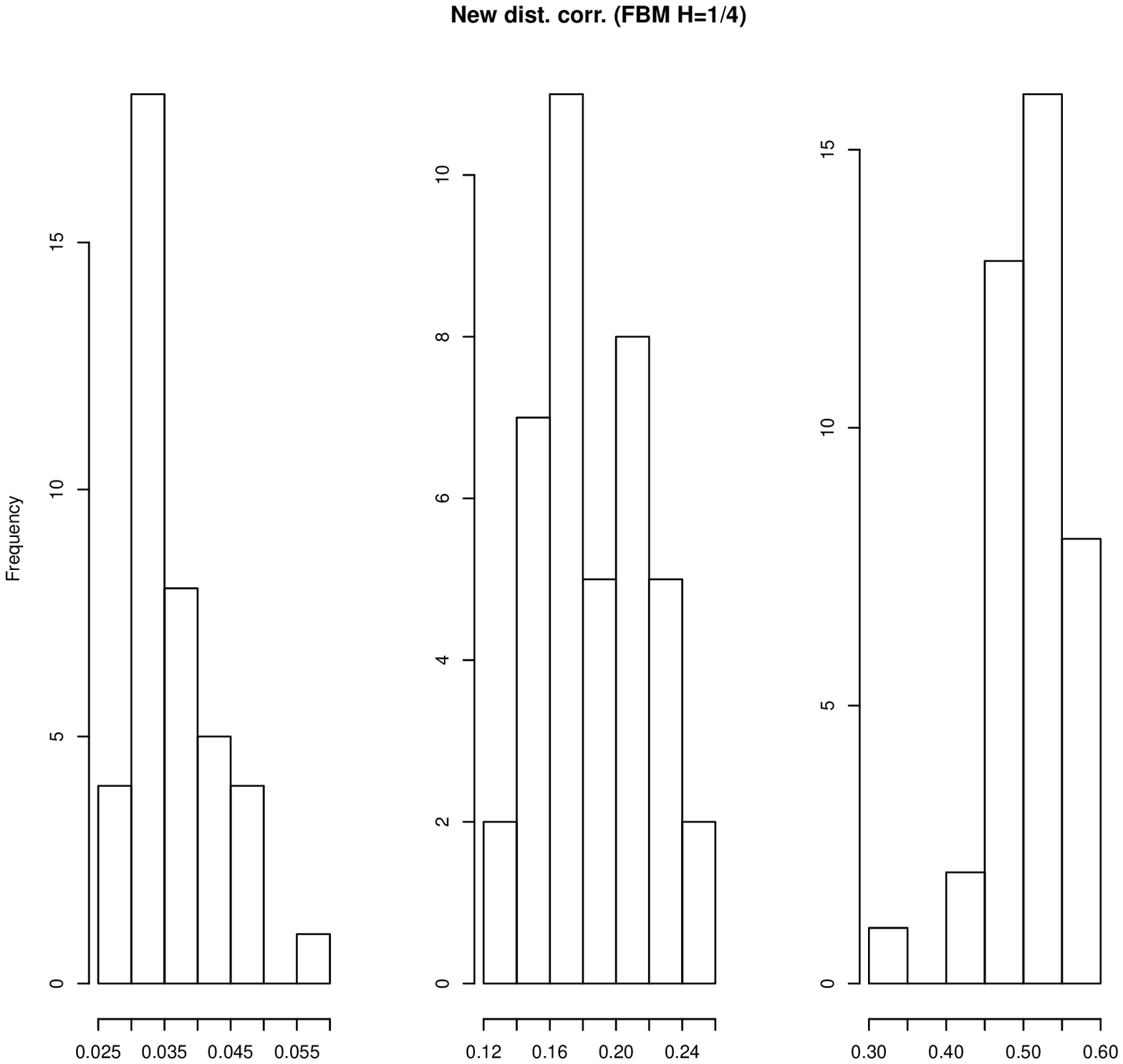}
\includegraphics[width=1\textwidth,height=10cm]{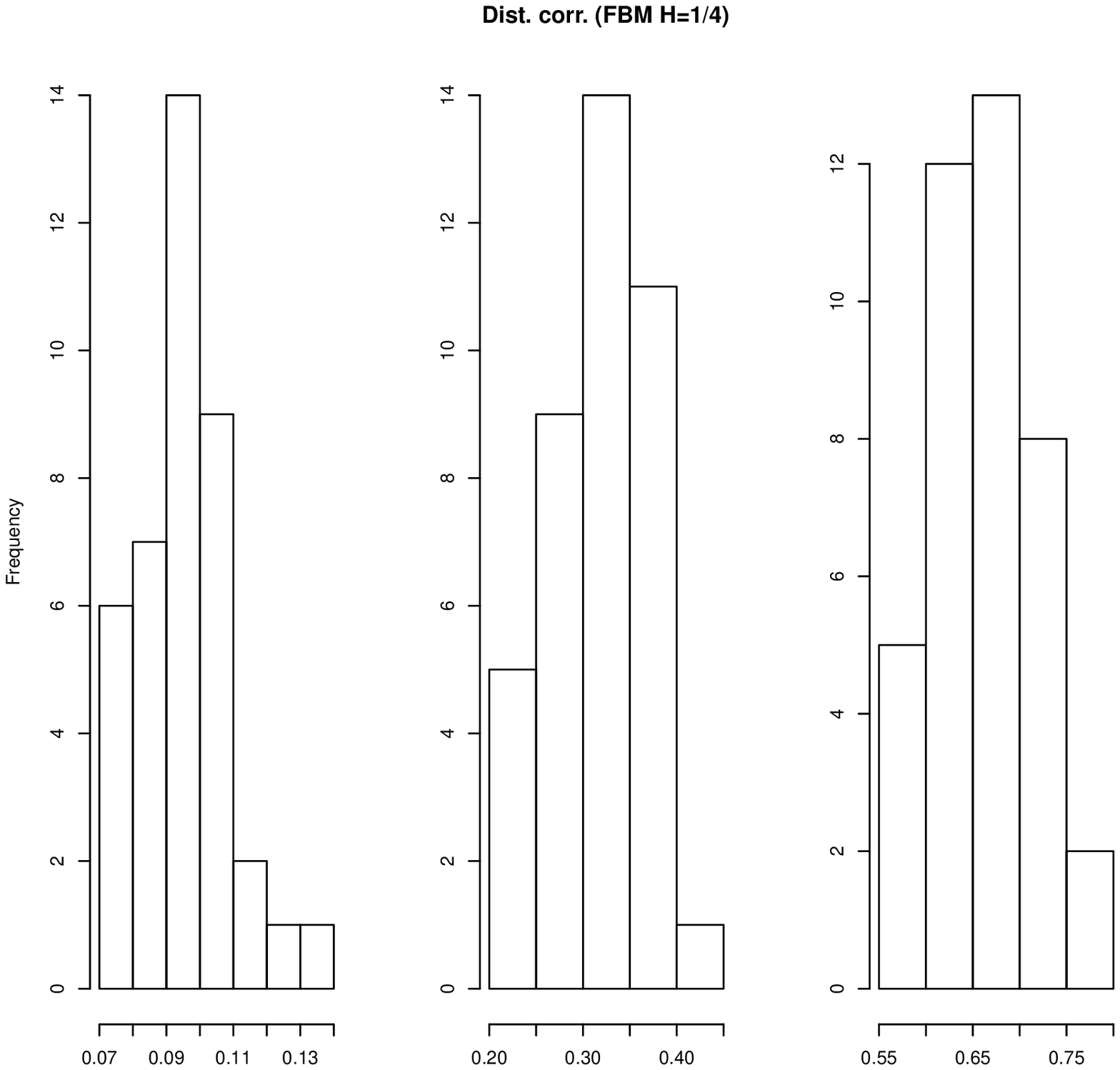} 
\caption{Histograms of $R_n(W_1,W_2)$ (top) and $R_n^{\rm Sz}(W_1,W_2)$ (bottom) for $H=0.25$ based on $40$ samples.  The 
 correlations of $W_1$ and $W_2$ are respectively $\rho=0,\,0.5,\,0.8$,
 from left to right.}
\label{fig:2}
\end{center}
\end{figure}
\begin{figure}[htbp]
\begin{center}
\includegraphics[width=1\textwidth,height=10cm]{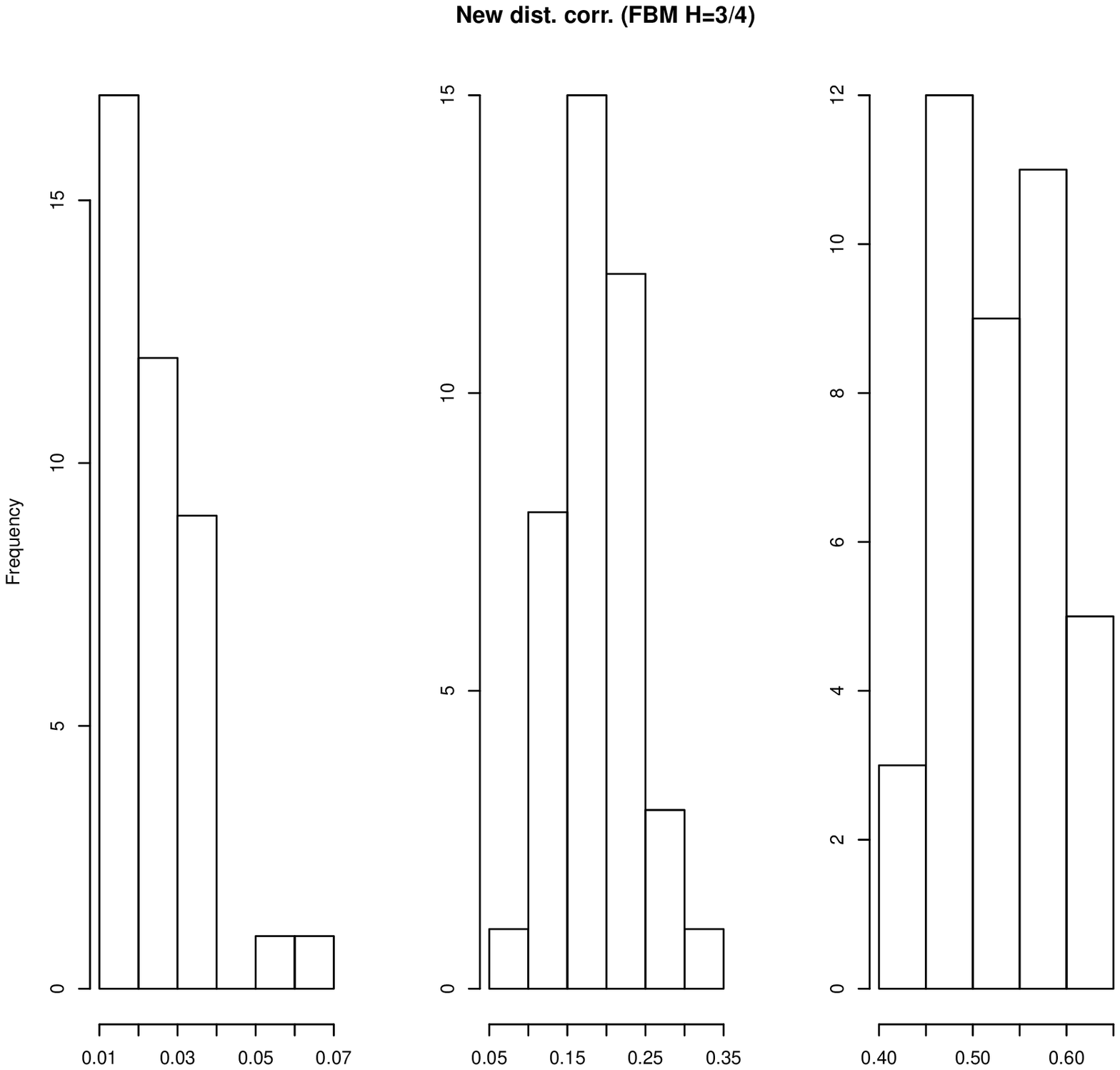}
\includegraphics[width=1\textwidth,height=10cm]{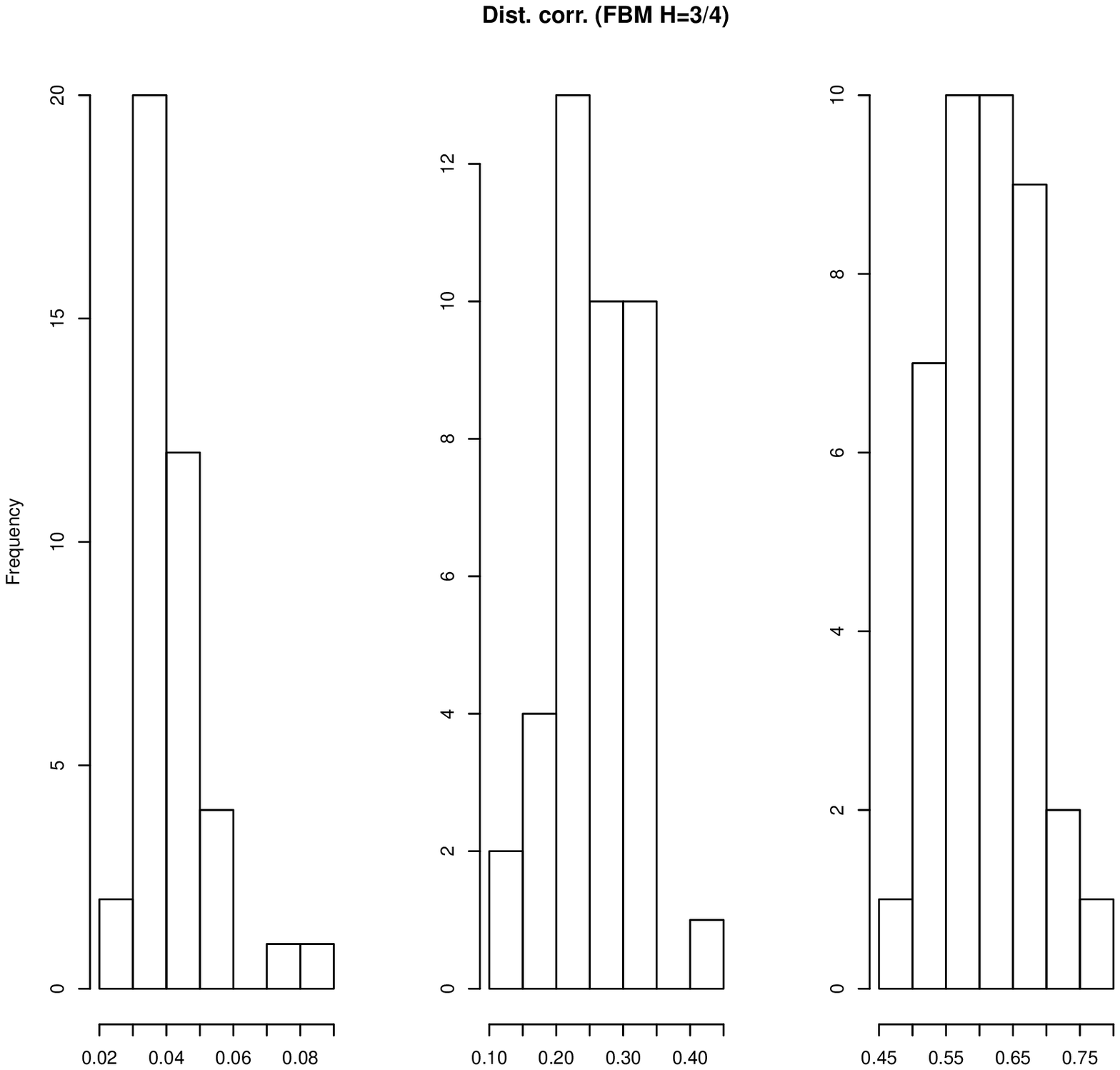} 
\caption{Histograms of $R_n(W_1,W_2)$ (top) and $R_n^{\rm Sz}(W_1,W_2)$ (bottom) for $H=0.75$ based on $40$ samples.  The 
correlations of $W_1$ and $W_2$ are respectively $\rho=0,\,0.5,\,0.8$,
from left to right.}
\label{fig:3}
\end{center}
\end{figure}

Figures \ref{fig:1}--\ref{fig:3} are based on 40 independent simulations of $R_n(X,Y)$ and $R_n^{\rm Sz}(\bfX,\bfY)$.
The 3 left (right) histograms show  $R_n(X,Y)$  ($R_n^{\rm Sz}(\bfX,\bfY)$) for 3 different choices of processes $(X,Y)$.
Although it is difficult to judge from such a small simulation study with rather special 
\spr es, these graphs give one the impression that both sample distance correlations capture the independence
or dependence of the processes $X$ and $Y$ quite well. The quantities $R_n^{\rm Sz}(\bfX,\bfY)$ have the tendency to
be larger than $R_n(X,Y)$.
\par
\begin{figure}[htbp] 
\begin{center}
\includegraphics[width=1\textwidth,height=7cm]{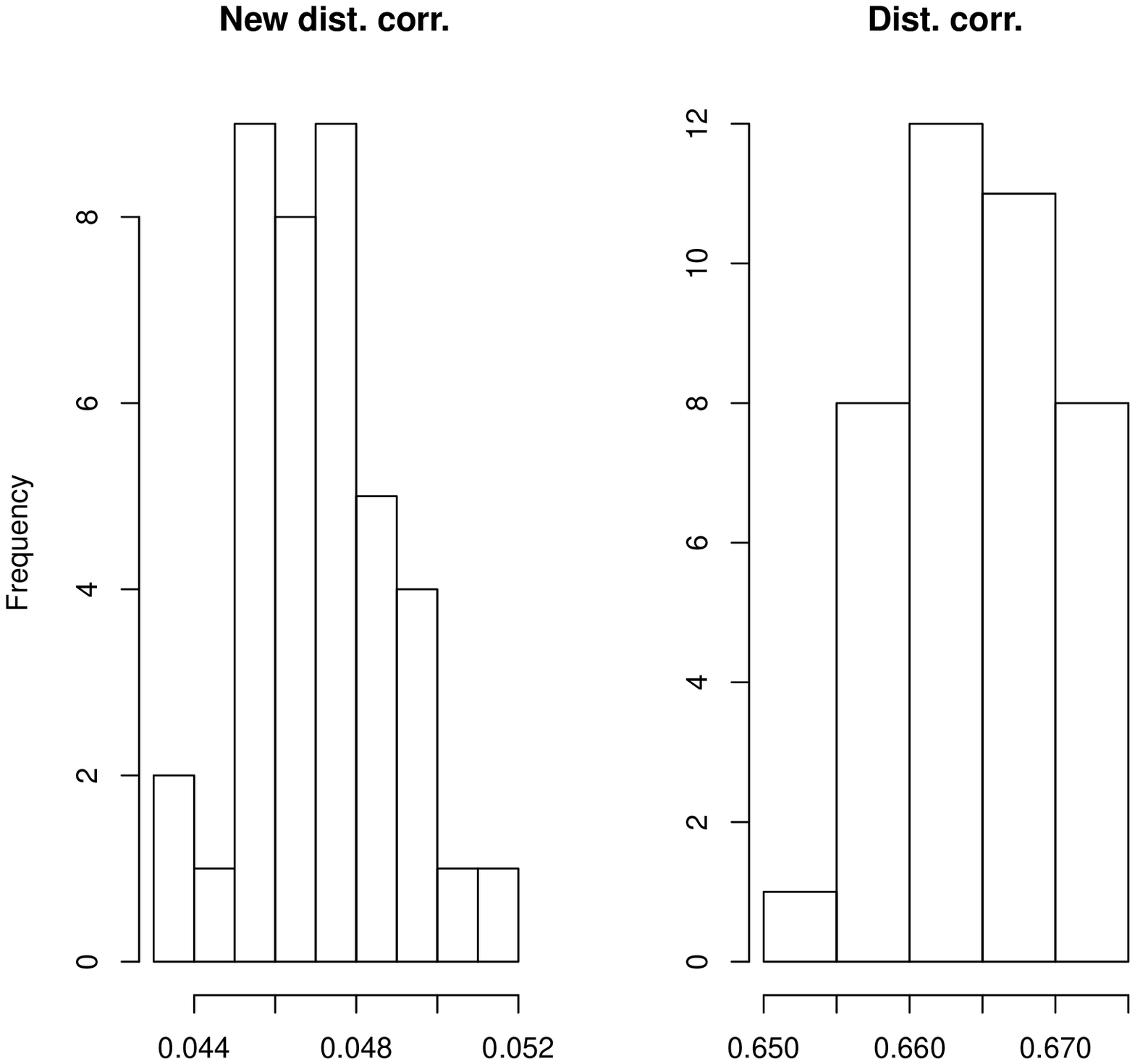}
\caption{Histograms of $R_n(X,Y)$ (left) and $R_n^{\rm Sz}(X,Y)$ (right) 
based on $40$ samples, where $X$ and $Y$ are independent piecewise
 constant processes based on iid normal \rv s. }
\label{fig:4}
\end{center}
\end{figure}

Finally, we consider two independent piecewise constant processes $X$ and $Y$ on $[0,1]$ which assume iid standard normal
values on the intervals $((i-1)/50,i/50],\,i=1,2,\ldots,50$. 
This is essentially the setting of
\cite{szekely:rizzo:2013} who chose independent vectors of iid normal random variables for the construction of
$R_n^{\rm Sz}(\bfX,\bfY).$ In the right histogram  of Figure \ref{fig:4} one can see that $R_n^{\rm Sz}(\bfX,\bfY)$ is typically far from zero.
This was observed in
\cite{szekely:rizzo:2013} who studied the case when the dimension of the vectors is large compared to the sample size.
On the other hand, our \ms\ $R_n(X,Y)$ is quite in agreement with the independence hypothesis.
\par
Of course, more investigations are needed in order to find out about the strengths and weaknesses of the 
distance covariances and correlation for processes introduced in this paper. One of the main
problems will be to find reliable confidence bands for the estimator $R_n(X,Y)$. This is work in progress.
\vspace{2mm} \\
\noindent {\bf Acknowledgment.} We would like to thank the referee for
constructive comments.

\end{document}